\documentclass[11pt]{article}
\font\tenmath=msbm10 scaled 1200
\font\sevenmath=msbm7 scaled 1200
\font\fivemath=msbm5 scaled 1200

\newfam\mathfam \textfont\mathfam=\tenmath
\scriptfont\mathfam=\sevenmath \scriptscriptfont\mathfam=\fivemath
\def\math{\fam\mathfam}

\def\N{{\math N}}
\def\E{{\math E}}

\def\P{{\math P}}

\def \^#1{\if#1i{\accent"5E\i}\else{\accent"5E#1}\fi}

\def \ind {1 \mkern -5mu \hbox{I}}

\def \a{\alpha}

\def \e{\varepsilon}
\def \g{\gamma}

\def \D{\Delta}

\def \cqfd{\quad_\diamondsuit} 
\def \ms{\medskip}
\def \ss{\smallskip}
\def \bs{\bigskip}
\def \ni{\noindent}
\def \F{{\cal F}}

\newtheorem{Lem}{Lemma}
\newtheorem{Pro}{Proposition}

\newtheorem{rem}{Remark}
\newenvironment{Rem}{\begin{rem}\rm}{\end{rem}}
\author{\sc Damien Lamberton \thanks{Laboratoire d'analyse et de math\'ematiques appliqu\'ees,
UMR~8050, Univ. Marne-la-Vall\'ee, Cit\'e
Descartes, 5 Bld Descartes, Champs-sur-Marne, F-77454 Marne-la-Vall\'ee Cedex 2.  {\tt
damien.lamberton@univ-mlv.fr}} 
\quad Gilles Pag\`es
\thanks{Laboratoire de probabilit\'es et mod\`eles  al\'eatoires,
UMR~7599, Univ. Paris~6, case~188, 4, pl. Jussieu, F-75252 Paris
Cedex 5. {\tt  gpa@ccr.jussieu.fr}} } 

\title{\bf How fast is the bandit?~\thanks{This work has benefitted from
the stay of both authors at the Isaac Newton Institute on the program {\em
Developments in Quantitative Finance}.}} 

\date{} 

\begin{document}

\maketitle
\begin{abstract}
In this paper we investigate the rate of convergence of the so-called two-armed bandit algorithm in a   financial context
of asset allocation. The behaviour of the algorithm turns out to be highly non-standard: no CLT whatever the time scale,
possible existence of two rate regimes.
\end{abstract}

\bigskip


\bs \ni {\em 2001 AMS classification:} 62L20,
secondary 93C40,
 91E40, 
68T05, 
 91B32 91B32
\section*{Introduction} 

In a recent joint work with P. Tarr\`es (see~\cite{LAPATA}), we studied the convergence of the so-called two-armed bandit
algorithm. In the terminology of learning theory (see $e.g.$~\cite{NA1, NA2}) this algorithm is a Linear Reward
Inaction ($LRI$) scheme. Viewed as a Markovian Stochastic Approximation ($SA$) recursive procedure, it   appears as the
simplest example of an algorithm having   two possible limits -- its {\em target} and a {\em
trap} -- both  {\em noiseless}. In
$SA$ theory   a target is a {\em stable} equilibrium of the Ordinary Differential Equation ($ODE$) associated to the {\em
mean} function of the algorithm, a trap being an unstable one. Various results from $SA$ theory show that an algorithm never
``falls" into a noisy trap (see~$e.g.$~\cite{LAZ, PEM, BRDU, FOPA,TAR}.    We established in~\cite{LAPATA} that
the two-armed bandit algorithm can be either {\em infallible} ($i.e.$ converging to its target with probability
one, starting from any initial value except the trap itself)  or {\em fallible}. This depends on the speed at
which the (deterministic)  {\em learning rate parameter} goes to
$0$.

Our aim on this paper is to investigate the rate of convergence of the algorithm,  toward either of its   limits. In fact,
the algorithm behaves in a highly non standard way  among $SA$ procedures. In particular, this rate is never ruled by a
Central Limit Theorem (CLT).    Furthermore, this study will provide some new insight on the infallibility problem as  it
will be seen further on.  However our motivations are not only theoretical but also practical in connection with the
financial context in which  the algorithm was presented in~\cite{LAPATA}, namely  a procedure for the optimal allocation of
a fund between the two traders who manage it. 
Imagine that the owner of a fund can share his wealth between two 
traders, say $A$ and $B$, and that, every day, he can evaluate the 
results of one of the traders and, subsequently,
modify 
the percentage of the fund managed by both traders. Denote by 
$X_n$ the percentage managed by trader $A$ at time $n$. We assume
that the owner selects the trader to be evaluated at random,
in such a way that the probability that $A$ is evaluated at time
$n$ is $X_n$, in order to select preferably the trader in charge of 
the greater part of the fund. In the $LRI$ scheme, if the evaluated 
trader performs well, its share is increased by a fraction
$\gamma_n\!\in(0,1)$ of the share of the other trader,  and nothing happens if
the evaluated trader performs badly. Therefore, the dynamics of the sequence
$(X_n)_{n\geq 0}$ can be modelled as follows:
\[
X_{n+1}=X_n+\gamma_{n+1}\left(
 \ind_{\{U_{n+1}\leq X_n\}\cap A_{n+1}}(1-X_n)
          -\ind_{\{U_{n+1}> X_n\}\cap B_{n+1}}X_n\right),\;\; X_0=x\!\in [0,1],
\]
where $(U_n)_{n\geq 1}$ is an i.i.d. sequence 
of uniform random variables on the interval $[0,1]$,
$A_n$ (resp. $B_n$) is the event ``trader $A$ (resp. trader $B$) performs well at 
time $n$". We assume $\P(A_n)=p_{\!_A}$, $\P(B_n)=p_{\!_B}$, for $n\geq 1$, 
with $p_{\!_A},\,p_{\!_B}\in(0,1)$, and independence between these events and the 
sequence $(U_n)_{n\geq 1}$. The point is that the owner of the fund does not know the parameters $p_{\!_A},\,p_{\!_B}$.
Note that this procedure is $[0,1]$-valued and that $0$ and $1$ are absorbing states. The
$\g_n$ parameter is the {\em learning} rate of the procedure (we will say from now on {\em reward} to take into account the
modelling context).

This recursive learning procedure   has been designed in order to assign  progressively 
  the whole fund to the best trader when $p_{_A}\neq p_{_B}$. From now on we will assume without loss of
generality that
$p_{_A}>p_{_B}$. This  means   that 
$X_n$ is expected to converge toward its target $1$ with probability $1$ provided
$X_0\!\in(0,1)$ (and consequently never to get trapped in 
$0$). However this ``infallibility" property needs  some very stringent
assumption on the reward  parameter $\g_n$: thus,  if $\g_n= \left(\frac{C}{C+n}\right)^\a$, $n\ge 1$, with $0<\a\le
1$ and
$C>0$, it is shown in~\cite{LAPATA} (see Corollary~1$(b)$) that  the algorithm is infallible    if and only if
$\a = 1$ and $C \le \frac{1}{p_{_B}}$.

In a standard $SA$ framework, when  an algorithm is converging to its target -- $i.e.$ a zero $x^*$
of its mean function $h(x) = \frac{\E (X_{n+1}-X_n\,|\,X_n=x)}{\g_n}$, stable for the $ODE$ $\dot x = h(x)$ --  its rate 
is ruled by a CLT at a $\sqrt{\g_n}$-rate with an asymptotic variance
$\sigma^2_{x^*}$ related to the asymptotic excitation of $x^*$ by the noise   
(see~\cite{BOU, KUYI, PEL}).  

As concerns the two-armed bandit algorithm, there is no exciting noise at $1$ (nor at $0$ indeed). This is made impossible
simply because both equilibrium points lie at the boundary of the state space $[0,1]$  of the algorithm (otherwise the
algorithm would leave the unit interval when getting too close to its boundary).  This same feature which    
causes the fallibility of the algorithm when $\g_n$ goes to $0$ too slowly also induces its non-standard rate of
convergence. 

To illustrate this behaviour and   consider again the steps $\g_n=
\frac{C}{C+n}$, $n\ge 1$,  with $C>0$. As a consequence of our main results, one obtains:
 
\begin{itemize}
\item If $C>\frac{1}{p_{_B}}$ the algorithm is fallible with positive probability
from any
$x\!\in[0,1)$ and, when failing, it goes to $0$ at a $n^{-Cp_{_B}}$-rate. The rate of convergence to $1$ may vary
according to the parameters, see Section~\ref{guidedtour}.  
\item If  $\frac{1}{p_{_A}-p_{_B}} \le  C\le \frac{1}{p_{_B}}$  (this case requires that $2\,p_{_B}\le p_{_A}$), the algorithm
is infallible from any $x\!\in (0,1]$ and goes to $1$ at a $n^{-Cp_{_A}}$-rate.
\item If  $\frac{1}{p_{_A}}< C<\frac{1}{p_{_A}-p_{_B}} $ then the algorithm is infallible (from any $x\!\in (0,1]$) and {\em
two rates of convergence to $1$ may occur with positive
$\P_x$-probability}: a ``slow" one -- $n^{-C(p_{_A}-p_{_B})}$~-- and a ``fast" one -- $n^{-Cp_{_A}}$.
\item If  $C\le \frac{1}{p_{_A}}$ then the algorithm is still infallible from any $x\!\in (0,1]$ but only
the slowest rate  of convergence ``survives" $i.e.$ $n^{-C(p_{_A}-p_{_B})}$.
\end{itemize}

In fact the following rule holds true: {\em   the greater the real constant $C$  is, the faster the algorithm
$(X_n)$ converges, except that when
$C$  is too great, then the algorithm becomes fallible} which makes the two-armed bandit a very ``moral"  procedure. 
Furthermore,  note that the ``blind" choice --
$C=1$ --  which ensures   infallibility    induces  a {\em slow rate} of convergence
$n^{-C(p_{_A}-p_{_B})}$ since then $C  \le
\frac{1}{p_{_A}}$ (by contrast with the {\em fast} rate $n^{-Cp_{_A}}$). 
 Also note that  this  rate is precisely that  of the {\em mean} algorithm
$x_{n+1}= x_n+\g_n(p_{_A}-p_{_B})x_n(1-x_n)$. A last feature to be noticed  is that the
switching between  rate regimes takes place ``progressively'' as the parameter $C$ grows since it happens that 
two different rates coexist  with positive probability. 

For more exhaustive results,  we refer to Section~\ref{guidedtour}. If one thinks again of a practical implementation of the
algorithm, the only reasonable choice for the reward parameter is $\g_n=\frac{1}{n+1}$: it ensures infallibility regardless
of the (unknown) values of
$p_{_A}$ and $p_{_B}$. But when these two parameters become too close, the rate of convergence becomes too poor to remain
really efficient. Unfortunately, this is more or less the standard situations: the daily performances of the traders are
usually close and this can be extended to other fields where this procedure can be used (experimental psychology, clinical
trials, industrial reliability,~\dots). One clue to get rid of this dependency is to introduce a   ``fading"
penalization in the procedure when an evaluated trader  has unsatisfactory performances. (By fading we mean
negligible with respect to the reward in order to preserve traders' motivation). This variant of the two-armed bandit
algorithm which satisfies a pseudo-CLT  at a (weak) $n^{-\frac 12}$-rate whatever the parameter $p_{_A}$ and
$p_{_B}$    is described and investigated in~\cite{LAPA5}. 

The paper is organized as follows: Section~\ref{one} is devoted to some preliminary results and technical tools.
Section~\ref{two} is devoted to the rate of convergence when the algorithm converges to its trap $0$ whereas 
Section~\ref{three} deals with the rate of convergence toward its target $1$. Section~\ref{guidedtour} proposes a summing up
of the results for a natural parameterized family of reward parameter $\g_n$.

\bs
\noindent {\sc Notations:}  
\noindent $\bullet$ Let $(a_n)_{n\ge 0}$ and $(b_n)_{n\ge 0}$ be two sequences of positive real numbers. The symbol
$a_n\sim b_n$ means $ a_n=b_n +o(b_n)$.

\smallskip
\noindent  $\bullet$  The notation $\P_x$ is used in reference to $X_0=x$.

 \section{Preliminary results}\label{one}
 We first recall the definition of the algorithm. We are interested 
 in the asymptotic behavior of the sequence $(X_n)_{n\in\N}$, where 
 $X_0=x$, with $x\in(0,1)$ and
 \[
 X_{n+1}=X_n+\gamma_{n+1}\left(
 \ind_{\{U_{n+1}\leq X_n\}\cap A_{n+1}}(1-X_n)
          -\ind_{\{U_{n+1}> X_n\}\cap B_{n+1}}X_n\right), \quad n\in \N.
\]
Here $(\gamma_n)_{n\geq 1}$ is a sequence of nonnegative numbers 
satisfying 
\[
\gamma_n<1\quad\mbox{  and }\quad  \Gamma_n=\sum_{k=1}^n \gamma_k\to +\infty\quad\mbox{  as }\quad n\to \infty,
\] $(U_n)_{n\geq 1}$ is 
a sequence of independent random variables which are uniformly 
distributed on the interval $[0,1]$, the events $A_n$, $B_n$
satisfy
\[
\P(A_n)=p_{_A},\quad \P(B_n)=p_{_B}, \quad n\in\N,
\]
where $0< p_{_B}<p_{_A}<1$, and the sequences $(U_n)_{n\geq 1}$
and $(\ind_{A_n}, \ind_{B_n})_{n\geq 1}$ are independent.
The natural filtration of the sequence $(U_n,\ind_{A_n}, \ind_{B_n})_{n\geq 1}$
is denoted by $(\F_n)_{n\geq 0}$ and
we set
\[
\pi=p_{_A}-p_{_B}>0.
\]
With this notation, we have, for $n\geq 0$,
\begin{equation}\label{mean-algo}
X_{n+1}=X_n+\gamma_{n+1}\pi X_n(1-X_n)+\gamma_{n+1}\Delta M_{n+1}, 
\end{equation}
where $\Delta M_{n+1}=M_{n+1}-M_n$, and the sequence $(M_n)_{n\geq 0}$ 
is the martingale defined by $M_0=0$ and
\[
\Delta M_{n+1}=\ind_{\{U_{n+1}\leq X_n\}\cap A_{n+1}}(1-X_n)
          -\ind_{\{U_{n+1}> X_n\}\cap B_{n+1}}X_n-\pi X_n(1-X_n).
\]
One derives from~(\ref{mean-algo}) that $(X_n)$ is a  $[0,1]$-valued super-martingale. Hence it converges $a.s.$ and
in $L^1$ to a limit $X_{_\infty}$. Consequently 
\[
\sum_n \g_nX_n(1-X_n) <+\infty \qquad a.s.
\]
which in turn shows that $X_{_\infty}=0$ or $1$ with probability $1$. One easily checks (see~\cite{LAPATA}) that $1$ is a 
stable equilibrium of the so-called {\em mean}  $ODE\equiv \dot x= \pi\,x(1-x)$ with attracting basin $(0,1]$ and $0$ is a
repulsive equilibrium of this $ODE$ (whence the terminology: $1$ is a {\rm target} and $0$ is a {\em trap}, see~\cite{LAPATA}
for more details).

\ss
The conditional variance process of the martingale $(M_n)$ will play
a crucial role in our analysis, and we will often use the following 
estimates.
\begin{Pro}\label{ProVarCond}
    We have, for $n\geq 0$,
    \[
    p_{_B} X_n(1-X_n)\leq \E\left(\D M_{n+1}^2\;|\;\F_n\right)
    \leq p_{_A} X_n(1-X_n).
    \]
\end{Pro}
\noindent {\sc Proof:} We have
\begin{eqnarray*}
     \E\left(\D M_{n+1}^2\;|\;\F_n\right)&=&
     p_{_A}X_n(1-X_n)^2+p_{_B}(1-X_n)X_n^2-\pi^2X_n^2(1-X_n)^2\\
     &=&X_n(1-X_n)\left(p_{_A}(1-X_n)+p_{_B}X_n-\pi^2X_n(1-X_n) \right)\\
     &\leq&X_n(1-X_n)\left(p_{_A}(1-X_n)+p_{_B}X_n \right)\\
     &\leq&p_{_A}X_n(1-X_n),
\end{eqnarray*}
where the last inequality follows from $p_{_B}\leq p_{_A}$. For the lower 
bound, note that
\begin{eqnarray*}
   p_{_A}(1-X_n)+p_{_B}X_n-\pi^2X_n(1-X_n)&=&(1-X_n)(p_{_A}-\pi^2X_n)+p_{_B}X_n\\
      &\geq &(1-X_n)(p_{_A}-\pi)+p_{_B}X_n=p_{_B},
\end{eqnarray*}
where we have used $\pi X_n\leq 1$.
$\cqfd$

\section{Convergence to the trap}\label{two}
We first prove that, under rather general conditions, as soon as the sequence converges to the
trapping state $0$, it goes to it very fast in the sense that the series $\sum_n X_n$ is convergent.

\begin{Pro}\label{EqFaillib}
If  
\begin{equation}\label{liminf-pi}
\liminf_n\frac{1}{\g_{n+1}}-\frac{1}{\g_n} > -\pi
\end{equation}
then
\[
\forall\, x\in (0,1),\qquad \{X_\infty=0\}=\{\sum_n X_n<+\infty\}\qquad \P_x\mbox{-}a.s. 
\]
\end{Pro}
Note that (\ref{liminf-pi}) is satisfied if the sequence 
$(\gamma_n)_{n\geq 1}$ is nonincreasing (for large enough $n$). 

\bs
\noindent {\sc Proof of Proposition \ref{EqFaillib}:}
Denote by $E$ the event $\{X_\infty=0\}\cap\{\sum_n X_n=+\infty\}$. We 
want to prove that $\P_x(E)=0$. We first show that on $E$,
\begin{equation}\label{mino-liminf}
\liminf_{n\to \infty}\frac{X_n}{\gamma_n\sum_{k=1}^nX_{k-1}}>0.
\end{equation}
We deduce from (\ref{mean-algo}) that 
\begin{eqnarray*}
    \frac{X_{n+1}}{\gamma_{n+1}}&=&\frac{X_n}{\gamma_{n+1}}+\pi 
    X_n(1-X_n)+\Delta M_{n+1}\\
    &=&\frac{X_n}{\gamma_n}+
    X_n\left(\frac{1}{\gamma_{n+1}}-\frac{1}{\gamma_n}
       +\pi (1-X_n)\right)+\D M_{n+1}.
\end{eqnarray*}
By summing up and setting $\gamma_0=\gamma_1$, we derive
\[
\frac{X_n}{\gamma_n}=\frac{x}{\gamma_1}+
   \sum_{k=1}^n \left(\frac{1}{\gamma_{k}}-
   \frac{1}{\gamma_{k-1}}
       +\pi (1-X_{k-1})\right)X_{k-1}+M_n.
       \]
 From Proposition~\ref{ProVarCond}, we know that the conditional 
 variance process of $(M_n)$ satisfies
 \[
 p_{_B}\sum_{k=1}^nX_{k-1}(1-X_{k-1})\leq \;<\!M\!>_n\;\leq p_{_A}\sum_{k=1}^nX_{k-1}(1-X_{k-1}).
 \]
 Therefore, on $E$, we have $<\!M\!>_\infty=+\infty$ a.s., and using
  the law of large numbers for martingales, we deduce that
 \[
 \lim_{n\to \infty}\frac{M_n}{\sum_{k=1}^nX_{k-1}}=0,\quad \mbox{a.s. on } E.
 \]
 The estimate (\ref{mino-liminf}) then follows easily from the 
 assumption~(\ref{liminf-pi}).
 
 Now let $S_n=\sum_{k=1}^n X_k$. Note that, on $E$, 
 $S_n\sim  \sum_{k=1}^n X_{k-1}$, so that, using~(\ref{mino-liminf}),
 \[
 \exists C>0,\quad\forall n\geq 1,\quad \gamma_n\leq C\frac{X_n}{S_n}.
 \]
 This implies
 \[
 \sum_n \gamma_n^2\leq C^2\sum_n \frac{X_n^2}{S_n^2}\leq 
 C^2\sum_n\frac{X_n}{S_n^2}<+\infty,
 \]
 where we have used $X_n\leq 1$. We also know from Proposition 9
 of \cite{LAPATA} (see (29) in particular)  that, on the set 
 $\{X_n\to 0\}$,
 \[
 \limsup_{n\to \infty}\frac{X_n}{\sum_{k\geq n}\gamma_{k+1}^2}<+\infty
 \quad\mbox{a.s.}
 \]
 Hence $X_n\leq C\sum_{k\geq n}\gamma_{k+1}^2$ for some $C>0$, and, 
 by plugging in the estimate $\gamma_{k+1}\leq CX_{k+1}/S_{k+1}$ we 
 derive
 \begin{eqnarray*}
     X_n&\leq &C\sum_{k\geq n}\frac{X_{k+1}^2}{S_{k+1}^2}\\
     &\leq &C\left(\sup_{k\geq n}X_{k+1}\right)
     \sum_{k\geq n}\frac{X_{k+1}}{S_{k+1}^2}\\
     &\leq&C\frac{\sup_{k\geq n}X_{k+1}}{S_n}.
     \end{eqnarray*}
     On the set $E$, we have $\displaystyle \lim_{n\to\infty}S_n=+\infty$, so, for 
     $n$ large enough, say $n\geq N$, we have 
     \[
     X_n\leq \frac{\sup_{k\geq n}X_{k+1}}{2}.
     \]
     Now, by taking $n$ to be the largest integer such that $X_n\geq 
     X_N$ (which exists on $\{X_n\to 0\}$ because $X_N>0$), we reach a 
     contradiction, which proves that $\P_x(E)=0$.
     $\cqfd$

\bigskip
Our next result shows that under~(\ref{ProConv0}), there is essentially only one way for $(X_n)$ to go to $0$.
\begin{Pro}\label{ProConv0}
    Assume~(\ref{liminf-pi}). 

\smallskip
\noindent $(a)$ Let $x\!\in(0,1)$. Then
\begin{equation}\label{Vers0ssi}
\P_x(X_\infty= 0)>0 \quad\Longleftrightarrow \quad\P(\sum_{n\ge 1} \prod_{k=1}^n(1-\mbox{\bf 1}_{B_k}\g_k)<+\infty)>0
\end{equation}
and,  on the event $\{X_\infty =0\}$, there exists a (random)  integer $n_0\ge 1$ such
that 
\begin{equation}\label{XVers0}
\forall\, n\ge n_0,\qquad X_n = X_{n_0}\prod_{k=n_0+1}^n(1-\mbox{\bf 1}_{B_k}\g_k).\qquad a.s.
\end{equation}
Note that, as a special case of~(\ref{Vers0ssi}),  
\begin{equation}\label{Vers0ssi2}
\sum_{n\ge 1}\prod_{k=1}^n(1-p_{_B}\g_k)<+\infty\quad
\Longrightarrow\quad \P_x(X_\infty = 0)>0.
\end{equation}
\smallskip
\noindent $(b)$ Furthermore, if $\displaystyle \;\sum_{n\ge 1 }\g^2_n <+\infty$, (\ref{Vers0ssi})   reads 
\[
\P_x(X_\infty = 0)>0\quad\Longleftrightarrow \quad \sum_{n\ge 1} \prod_{k=1}^n(1-p_{_B}\g_k)<+\infty 
\]
 and moreover there is a  random
variable $ \Xi_x>0$ such that 
\[
X_n\sim \Xi_x\prod_{k=1}^n(1-p_{_B}\g_k)\quad \qquad a.s.\quad \mbox{ on }\quad 
\{X_\infty = 0\}. 
\]
\end{Pro}

\begin{Rem} If  $\displaystyle \sum_{n\ge 1 }\g^2_n =+\infty$,  a weaker (but still tractable)
sufficient condition for $\P_x(X_\infty = 0)$ is given by 
\[
\exists \, \rho\!\in(0,p_{_B}(1-p_{_B})/2),\qquad \sum_{n\ge 1}e^{-\rho\Gamma^{(2)}_n}
\prod_{k=1}^n(1-p_{_B}\g_k)<+\infty
\] 
where $\Gamma^{(2)}_n = \sum_{1\le k\le n}\g^2_k$ (see the proof of Proposition~\ref{ProConv0}). Then, on the set
$\{X_n\to 0\}$, for every $\eta\!\in(0,p_{_B}(1-p_{_B})/2)$, 
\[
X_n = o\left(e^{-(\frac{p_{_B}(1-p_{_B})}{2}-\eta)\Gamma^{(2)}_n}\prod_{k=1}^n(1-p_{_B}\g_k)\right).
\]
\end{Rem}
\begin{Rem}
 Note that the condition in~(\ref{Vers0ssi}) which characterizes 
fallibility does not depend on
$x$: if the algorithm is fallible for one $x\!\in(0,1)$ then it is for any such $x$. 
\end{Rem}
\bigskip
\noindent {\sc Proof of Proposition~\ref{ProConv0}:}
$(a)$ It follows from Proposition~\ref{EqFaillib} and the conditional Borel-Cantelli Lemma that
$\P_x\mbox{-}a.s.$
\begin{equation}\label{Egalite1}
\{X_n\to 0\}= \{\sum_{n\ge 0} \mbox{\bf 1}_{\{U_{n+1}\le X_n\}}<+\infty\}=
\bigcup_{n\ge 0}\bigcap_{k\ge n} \left\{U_{k+1} >X_{k}\right\}.
\end{equation} 
The sequence of events $\left(\bigcap_{k\ge n} \left\{U_{k+1} 
>X_{k}\right\}\right)_{n\geq 1}$ being non-decreasing, we have
\[
\P_x(X_n\to 0)=\lim_{n\to \infty}\P_x\left(\bigcap_{k\ge n} 
\left\{U_{k+1} >X_{k}\right\}\right),
\]
and the left-hand side is positive if and only if, for some integer 
$n\geq 1$,
\[
\P_x\left(\bigcap_{k\ge n} 
\left\{U_{k+1} >X_{k}\right\}\right)>0.
\]
From the definition of the sequence $(X_n)$, we get
 (with the convention $\prod_{\emptyset}=1$), 
\begin{eqnarray} 
\bigcap_{k\ge n} 
\left\{U_{k+1} >X_{k}\right\}&=& \bigcap_{k\ge
n}\left\{U_{k+1}>X_k \mbox{ and }X_{k}=X_n\prod_{\ell=n+1}^k\hskip -0.25 
cm (1-\mbox{\bf 1}_{B_\ell}\g_\ell)\right\}\label{Descente}\\
&=&\bigcap_{k\ge n}\left\{U_{k+1}>X_n\hskip -0.25 cm \prod_{\ell=n+1}^k(1-\mbox{\bf
1}_{B_\ell}\g_\ell)\right\}.
\label{Inclu2} \end{eqnarray}
Note that (\ref{XVers0}) follows from (\ref{Egalite1}) and (\ref{Descente}).
Now, denote by ${\cal B}_n$ the $\sigma$-field generated by the random 
variable $X_n$ and the events, $B_k$, $k\geq n$. We have
\[
\P_x\left(\bigcap_{k\ge n}\left\{U_{k+1}>X_n\hskip -0.25 cm 
\prod_{\ell=n+1}^k(1-\mbox{\bf
1}_{B_\ell}\g_\ell)\right\}\;|\; {\cal B}_n\right)=
\prod_{k=n}^\infty\left(
        1-X_n\prod_{l=n+1}^k(1-\mbox{\bf
1}_{B_\ell}\g_\ell)\right),
\]
and the infinite product is positive if and only if
\[
\sum_k\prod_{l=n+1}^k(1-\mbox{\bf
1}_{B_\ell}\g_\ell)<+\infty.
\]
This clearly implies (\ref{Vers0ssi}).
 The  sufficient condition~(\ref{Vers0ssi2}) follows from the equality 
\[
\E\left(\sum_{n\ge 1}\prod_{1\le k\le n}(1-\mbox{\bf 
1}_{_{B_k}}\g_{k})\right)= \sum_{n\ge 1}\prod_{1\le k\le n}(1-p_{_{B}}\g_{k}).
\]

\smallskip
\noindent $(b)$ (and proof of the remark) If $\displaystyle \sum_{n\ge 1}\g^2_n<+\infty$, then,
a straightforward argument (see~\cite{LAPATA}, proof of Lemma~2) shows 
that
\[
\prod_{k=1}^n \left(\frac{1-\mbox{\bf 1}_{_{B_k}}\g_{k}}{1-p_{_{B}}\g_{k}}\right)\longrightarrow \xi>0
\qquad\mbox{a.s.}\qquad n\to+\infty.
\]
This proves   claim~$(b)$. 

When $\displaystyle \sum_{n\ge 1}\g^2_n=+\infty$, 
one checks that 
\[
\log \prod_{k=1}^n \left(\frac{1-\mbox{\bf 1}_{_{B_k}}\g_{k}}{1-p_{_{B}}\g_{k}}\right)= M^B_n-\sum_{k=1}^n(\frac 12
p_{_B}(1-p_{_B})+\varepsilon_k)\g^2_k .
\]
where $\varepsilon_k$ is random variable bounded by $c\g_k$ ($c$ real constant) and 
\[
M^B_n = \sum_{k= 1}^n(\mbox{\bf 1}_{_{B_k}}-p_{_B})\g_k(1-\g_k/2)
\]
is a martingale with bounded increments satisfying $<\!M^B\!>_{n}\sim p_{_B}(1-p_{_B})  \Gamma^{(2)}_n\to+\infty$. Then
\[
M^B_n = o\left(\Gamma^{(2)}_n\right)
\]
since $\frac{M^B_n}{<\!M^B\!>_n}\to 0$ as $n\to\infty$. Consequently, $\P$-$a.s.$, there exists a finite random variable $\xi$ such
that 
\[
\prod_{k=1}^n  (1-\mbox{\bf 1}_{_{B_k}}\g_{k}) \le \xi \exp{\left(-(\frac 12
p_{_B}(1-p_{_B})+o(1))\Gamma^{(2)}_n\right)}\prod_{k=1}^n (1-p_{_{B}}\g_{k})
\]
where $o(1)$ denotes a random variable $\P$-$a.s.$  going to $0$ as $n\to \infty$. The sufficient condition given in the
remark follows straightforwardly as well as  the rate of convergence of $X_n$.$\cqfd$

\section{Convergence to the target}\label{three}
In order to study the rate of convergence to $1$, we first rewrite 
(\ref{mean-algo}) as follows:
\begin{equation}\label{mean-algo-to1}
1-X_{n+1}=(1-X_n)\left(1-\gamma_{n+1}\pi X_n\right)-\gamma_{n+1}\D M_{n+1}.
\end{equation}
Now let
\[
\theta_n=\prod_{k=1}^n (1-\gamma_k\pi X_{k-1}), \quad 
Y_n=(1-X_n)/\theta_n, \quad n\in\N.
\]
\begin{Pro}
    \label{Fact1} $(a)$
    The sequence $(Y_n)_{n\in\N}$ is a non-negative martingale.
     
\ss
\noindent $(b)$ 
        On the set $\{X_\infty=1\}$, we have 
        \[\lim_{n\to \infty}\frac{1-X_n}{\prod_{k=1}^n(1-\pi\gamma_k)}=
          \xi Y_\infty
        \]
        almost surely, where $\xi$ is a finite positive random 
        variable and $Y_\infty=\lim_{n\to\infty}Y_n$.
 
\end{Pro}
\noindent {\sc Proof:} The first assertion follows from the equality
\[
Y_{n+1}=Y_n-\frac{\gamma_{n+1}}{\theta_{n+1}}\D M_{n+1},
\]
and the fact that the sequence $(\theta_n)_{n\in\N}$ is predictable.

As a non-negative martingale, the sequence $(Y_n)_{n\in\N}$ has a limit
$Y_\infty$, which satisfies $Y_\infty\geq 0$ a.s. and 
$\E(Y_\infty)<+\infty$. 

Recall that $\sum_n \gamma_n  X_{n-1}(1-X_{n-1})<+\infty$ almost 
surely. Therefore, on $\{X_\infty=1\}$, we have $\sum_n \gamma_n  
(1-X_{n-1})<+\infty$ a.s.,
which implies that the sequence 
$\prod_{k=1}^n\frac{1-\pi \gamma_k X_{k-1}}{1-\pi \gamma_k}$ has a 
positive and finite limit
and the second assertion of the Proposition 
follows easily.$\cqfd$

\begin{Rem}\label{rem-theta}
    Note that, with the notation $\Gamma_n=\sum_{k=1}^n\gamma_k$,
    we have 
    $\prod_{k=1}^n(1-\pi \gamma_k)\leq e^{-\pi \Gamma_n}$.
    Therefore, we deduce from Proposition~\ref{Fact1} that, on the 
    set $\{X_\infty=1\}$, $1-X_n=O(e^{-\pi \Gamma_n})$ almost surely.
    If we have 
    $\sum_n\gamma_n^2<+\infty$, the sequence 
    $\left( e^{\pi \Gamma_n}\prod_{k=1}^n(1-\pi \gamma_k)\right)$
     converges to a 
    positive limit, so that, on the set $\{X_\infty=1\}$, 
    we have $\displaystyle \lim_{n\to \infty}e^{\pi \Gamma_n}(1-X_n)=\xi'Y_\infty$,
    with $\xi'\in (0,+\infty)$ almost surely. 
    
    On the other hand, on $\{X_\infty=0\}$,
    the sequence $(\theta_n)_{n\in\N}$ itself converges to an almost 
    surely positive limit, so that 
    $\{Y_\infty=0\}\subset \{X_\infty=1\}$.
    \end{Rem}
 \begin{Pro}\label{Fact2}$(a)$  If $\sum_n \gamma_n^2e^{\pi \Gamma_n}<+\infty$, the martingale
     $(Y_n)_{n\in\N}$ is bounded in $L^2$ and its limit satisfies
     $\E(X_\infty Y_\infty)>0$. Moreover, on the set $\{Y_\infty=0\}$,
     we have 
     \begin{equation}\label{Ysuretc}
     \limsup_{n\to\infty}\frac{Y_n}{
     \sum_{k\geq n}\gamma_{k+1}^2e^{\pi\Gamma_{k+1}}}<+\infty
     \end{equation}
     almost  surely.
     
\ss
\noindent $(b)$   If  
\begin{equation}\label{Xvers1a}
\sum_n\g^2_n\,e^{\pi \Gamma_n}=+\infty\quad\mbox{and}\quad
\sup_{n\geq 1}\gamma_ne^{\pi\Gamma_n}<+\infty,
\end{equation}
then, for every $x\!\in(0,1)$, 
\[
\{X_\infty = 1\} =\{Y_{_\infty}=0\}\quad \P_x\mbox{-}a.s.
\]
 \end{Pro}

 \begin{Rem}
 It follows from Proposition \ref{Fact2} and Remark \ref{rem-theta}
 that, if $\sum_n \gamma_n^2e^{\pi \Gamma_n}<+\infty$, on the set
 $\{X_\infty=1\}\cap\{Y_\infty>0\}$ (which has positive probability)
 the sequence $((1-X_n)e^{\pi \Gamma_n})_{n\in\N}$ converges to a 
 positive limit almost surely. 
 \end{Rem}
 \begin{Rem}
 We also derive from the inequality
 $(1-X_{n+1})\geq (1-X_n)\left(1-\gamma_{n+1}\ind_{\{U_{n+1}\leq X_n\}\cap 
 A_{n+1}}\right)$ that 
$$
1-X_n\geq 
 (1-x)\prod_{k=1}^n(1-\gamma_k\ind_{A_k})\geq Ce^{-p_{_A}\Gamma_n},
$$ 
for 
 some real constant $C>0$, if
 $\sum_n\gamma_n^2<+\infty$. Therefore, we deduce from 
 Proposition~\ref{Fact2} that if 
     $\displaystyle \lim_{n\to \infty}\left(e^{p_{_B}\Gamma_n}
     \sum_{k\geq n}\gamma_{k+1}^2e^{\pi\Gamma_{k+1}}\right)=0$, then 
 $\P(Y_\infty=0)=0$. On the other hand, the second part of 
 Proposition~\ref{Fact2} shows that, in some cases, we may
 have $1-X_n=o(e^{-\pi\Gamma_n})$, and we need to investigate what
 the real rate of convergence is in such cases: see 
 Proposition~\ref{Fact3}.
 \end{Rem}
 
 \ms
 
 \noindent{\sc Proof of Proposition~\ref{Fact2}:} $(a)$ Assume $\sum_n \gamma_n^2e^{\pi \Gamma_n}<+\infty$.
 In order to prove $L_2$-boundedness, we estimate the conditional
 variance process. Using Proposition~\ref{ProVarCond}, we have
 \begin{eqnarray}
 	\nonumber    \E\left( (Y_{n+1}-Y_n)^2 \;|\; 
     \F_n\right)&=&\frac{\g_{n+1}^2}{\theta_{n+1}^2}
          \E\left(\D M_{n+1}^2\;|\;\F_n\right)\\
		\nonumber&\leq &\frac{\g_{n+1}^2}{\theta_{n+1}^2}p_{_A}X_n(1-X_n)\\
	\nonumber	&= &\frac{\g_{n+1}^2}{\theta_{n+1}^2}p_{_A}X_n\theta_nY_n\\
	\nonumber	&\leq &p_{_A}\frac{\g_{n+1}^2}{\theta_{n}(1-\pi \g_{n+1})^2}\,Y_n\\
	\nonumber &\leq &p_{_A}\frac{\g_{n+1}^2}{(1-\pi \g_{n+1})^2\prod_{k=1}^n(1-\pi 
	\g_{k})}\,Y_n\\
	\label{crochetY}&\leq & C\,p_{_A} \g_{n+1}^2 e^{\pi\Gamma_{n+1}}Y_n,
 \end{eqnarray}
 where we have used the inequality $\theta_{n}\geq \prod_{k=1}^n(1-\pi \g_{k})$
 and the fact that, since we have $\sum_{n\geq 1}\g_n^2<+\infty$,
 $\prod_{k=1}^n (1-\pi \g_{k})\geq  e^{-\pi\Gamma_n}/C$ for some  $C>0$.
 Note that $\displaystyle \sup_{n\in\N}\E \,Y_n<+\infty$. Therefore, the convergence
 of the series $\sum_n\g_n^2 e^{\pi\Gamma_{n}}$ implies that 
 $(Y_n)_{n\in\N}$ is bounded in $L_2$.
 
 In order to prove $\E(X_\infty Y_\infty)>0$, we consider the 
 conditional covariance
 \begin{eqnarray*}
\E_x\,((1-X_n)X_n\,|\, {\cal F}_{n-1})&=& X_{n-1}(1-X_{n-1})
\left(1+\pi\, \g_n(1-2X_{n-1})+\pi \,\g_n^2X_{n-1}-p_{_A}\g^2_n\right)\\
&\ge& X_{n-1}(1-X_{n-1})\left(1-\pi\,\g_n X_{n-1}-p_{_A}\g^2_n\right)
\end{eqnarray*}
\begin{eqnarray*}
\mbox{so that }\hskip 2,5 cm  \E_x\left(\,X_nY_n\,|\, {\cal F}_{n-1}\right)&\ge&
X_{n-1}Y_{n-1}\left(1-\frac{p_{_A}\g_n^2}{1-\pi
\g_n  X_{n-1}}\right) \hskip 2,5 cm\\
     &\geq &X_{n-1}Y_{n-1}\left(1-\frac{p_{_A}\g_n^2}{1-\pi \g_n}\right).
\end{eqnarray*}
For $n$ large enough (say $n\geq n_0$), we have 
$1>\frac{p_{_A}\g_n^2}{1-\pi \g_n}$ and,
 by  induction, for $n\geq n_0$,
\begin{eqnarray*}
\E_x(X_nY_n)\ge \E_x  X_{n_0}Y_{n_0} \prod_{k=n_0+1}^n\left(1-\frac{p_{_A}\g_k^2}{1-\pi \g_k}\right).
\end{eqnarray*}
Now, using that $Y_n\to Y_{_\infty}$ and $X_n\to X_{_\infty}$ in $L^2(\P)$,
and $\sum_n\g_n^2<+\infty$, one finally gets
$\E_x(X_{_\infty}Y_{_\infty})>0$.
Note that
 this implies that $\P_x(X_{_\infty}=1,\, Y_{_\infty}>0)>0$ since
$X_{_\infty}=\mbox{\bf 1}_{\{X_{_\infty}=1\}}$.

 The first step to establish~(\ref{Ysuretc}) is  to apply to the martingale
$(Y_n)_{n\ge 1}$ an approach originally developed in~\cite{LAPATA} to establish the infallibility property for
$(X_n)$:  for every $n\ge 1$, 
\begin{eqnarray*}
\P(Y_{_\infty}=0\,|\,{\cal F}_n)&= &\frac{1}{Y_n^2}\,\E_x(\mbox{\bf 1}_{\{Y_{_\infty}=0\}}(Y_{_\infty}-Y_n)^2\,|\,{\cal F}_n) \\
&\le& \frac{1}{Y_n^2}\sum_{k\ge n+1}\E_x( (Y_{k}-Y_{k-1})^2\,|\,{\cal F}_n).
\end{eqnarray*}
Plugging~(\ref{crochetY}) in the above inequality  and using that $\E_x(Y_k\,|\,{\cal F}_n) =Y_n$ for every $k\ge n$ yield, 
\[
\P_x(Y_{_\infty}=0\,|\,{\cal F}_n)\le \frac{C p_{_A}}{Y_n}\sum_{k\ge n+1}\g^2_k e^{\pi \Gamma_k}.
\]
On the other hand the martingale $\P_x(Y_{_\infty}=0\,|\,{\cal F}_n)$  converges $\P_x$-$a.s.$ toward $\mbox{\bf
1}_{\{Y_{_\infty}=0\}}$. The announced result follows easily.

\ss
\noindent 
$(b)$  We now assume $\sum_n \gamma_n^2e^{\pi \Gamma_n}=+\infty$ and
$\displaystyle \sup_n \gamma_ne^{\pi\Gamma_n}<+\infty$. Note that the latter 
condition implies $\g_n^2\leq C\gamma_ne^{-\pi\Gamma_n}$ for some 
$C>0$,
so that $\sum_n\g_n^2<+\infty$. On the other hand, we have
\begin{eqnarray*}
    |Y_n-Y_{n-1}|&=&\frac{\g_n}{\theta_n}|\D M_n|\\
           &\leq&\frac{\g_n}{\prod_{k=1}^n(1-\pi\g_k)}|\D M_n|\leq 
           C\g_ne^{\pi \Gamma_n}|\D M_n|,
 \end{eqnarray*}
 so that the martingale $(Y_n)_{n\geq 1}$ has bounded increments.
 Consequently the Law of
Iterated Logarithm (cf.~\cite{HAHE}) implies that $\displaystyle \liminf_n Y_n =-\infty$ on the event
$\{<\!Y\!>_{_\infty}=+\infty\}$, and, since $Y_n\geq 0$, we deduce 
thereof that $\{<\!Y\!>_{_\infty}<+\infty\}$ almost surely. On the other hand, 
we have, using Proposition~\ref{ProVarCond} and the inequality 
$\theta_n\leq e^{-\pi\Gamma_n}$, 
\begin{eqnarray*}
   \Delta\!<\!Y\!>_n&=&\frac{\g_n^2}{\theta_n^2}
   \E\left(\D M_n^2\;|\;\F_{n-1}\right)\\
        &\geq &\frac{\g_n^2}{\theta_n(1-\pi\gamma_n)}p_{_B}X_{n-1}Y_{n-1}\\
       &\geq&CX_{n-1}Y_{n-1}\g_n^2e^{\pi \Gamma_n}.
\end{eqnarray*}   
Therefore,  the assumption (\ref{Xvers1a}) implies that $Y_{\infty}=0$ 
on the event $\{X_\infty=1\}$.\hfill$\cqfd$

\ms

In order to clarify what happens when $Y_\infty=0$, we first observe that
we have, up to null events,
\begin{eqnarray*}
\left\{\sum_n(1-X_n)<+\infty\right\}&=&\left\{\sum_n\ind_{\{U_n>X_n\}}<+\infty\right\}\\
     &\subset&
     \bigcup_{m\geq 1}\bigcap_{n\geq 
     m}\left\{1-X_n=(1-X_{m})\prod_{k=m+1}^n(1-\ind_{A_k}\g_k)\right\},
\end{eqnarray*}
so that, on the set $\left\{\sum_n(1-X_n)<+\infty\right\}$, we have
\[
1-X_n\sim \xi\prod_{k=1}^n(1-\ind_{A_k}\g_k)\quad\mbox{a.s.},
\]
where $\xi$ is a positive random variable. Recall that, if 
$\sum_n\g_n^2<+\infty$, $\prod_{k=1}^n(1-\ind_{A_k}\g_k)\sim 
\xi'e^{-p_{_A}\Gamma_n}$, for some (random) $\xi'>0$.
We thus see that, on the set 
$\left\{\sum_n(1-X_n)<+\infty\right\}$, we have a ``fast" rate 
of convergence. 
The possibility of occurrence of this fast rate is characterized
in the following Proposition.
\begin{Pro}\label{vitevitevers1}
    We have, for all $x\in(0,1)$,
    \[
    \P_x(\sum_n(1-X_n)<+\infty)>0
    \quad\Longleftrightarrow \quad\P(\sum_{n\ge 1} 
    \prod_{k=1}^n(1-\mbox{\bf 1}_{A_k}\g_k)<+\infty)>0.
    \]
\end{Pro}
Note that the condition $\displaystyle\sum_{n\ge 1} \prod_{k=1}^n(1-p_{_A}\g_k)<+\infty$
implies $\displaystyle\P(\sum_{n\ge 1} 
    \prod_{k=1}^n(1-\mbox{\bf 1}_{A_k}\g_k)<+\infty)=1$ and that if 
    $\sum_n\gamma_n^2<+\infty$, we have
    \[
    \P_x(\sum_n(1-X_n)<+\infty)>0
    \quad\Longleftrightarrow \quad\sum_{n\ge 1}e^{-p_{_A}\Gamma_n}<+\infty.
    \]
    The proof of Proposition~\ref{vitevitevers1} and of these comments 
    is similar to that of the analogous statements concerning 
    convergence to 0.

    \ms
    
    In the following Proposition, we give a sufficient condition
    for the fast rate to be achieved with probability one and 
    a sufficient condition under which we have at most two rates with 
    positive probability: $e^{-\pi\Gamma_n}$ and 
    the fast rate $e^{-p_{_A}\Gamma_n}$.
 \begin{Pro}\label{Fact3}
     Let $\e_n=\frac{1}{\gamma_{n+1}}-\frac{1}{\gamma_n}-\pi$ for 
					$n\ge 1$. 
     
\ss
\ni $(a)$ If $\sum_n\gamma_n\e_n^+<+\infty$, we have
     $\sum_n(1-X_n)<+\infty$ almost surely on the set $\{X_\infty=1\}$.
    
\ms
\ni $(b)$ 
     If $\displaystyle \liminf_n \e_n>0$, then $\sum_n\gamma_n^2e^{\pi 
     \Gamma_n}<+\infty$, and, on the event $\{Y_\infty=0\}$,
     we have  $\sum_n(1-X_n)<+\infty$ almost surely.
     
 \end{Pro}
 Note that the condition $\sum_n\gamma_n\e_n^+<+\infty$ implies $\displaystyle \liminf_n \varepsilon_n^+=0$ and is satisfied in 
 the following
 cases:
 \begin{itemize}
     \item the sequence $(\gamma_n)$ is constant,
     \item $\gamma_n=\lambda n^{-\alpha}$ (for large enough $n$), with $\lambda$ a positive 
     constant and $0<\alpha<1$,
     \item $\gamma_n=C/(C+n)$, where the constant $C$ satisfies      
     $\pi C\geq 1$.
 \end{itemize}
 On the other hand, if $\gamma_n=C/(C+n)$, with $\pi C<1$, we have
  $\displaystyle \liminf_n \e_n>0$.
  
  \medskip
 
  Before proving Proposition~\ref{Fact3}, we state and prove a lemma which will
  be useful for the proof of the second statement.
  \begin{Lem}\label{LemmeFast}
      Assume that, for some positive  integer $n_0$, 
      $
      \forall n\geq n_0,\quad \e_n\geq 0.
      $
      Then, the sequence $(Z_n)_{n\geq n_0}$, with 
      $Z_n=(1-X_n)/\g_n$ is a submartingale, and we have $\sum_n(1-X_n)<+\infty$ a.s., on the set 
      $\{X_\infty=1\}\cap\{\sup_n\frac{1-X_n}{\g_{n+1}}<+\infty\}$.
  \end{Lem}
  \begin{Rem}
      If $\displaystyle \inf_n\g_n e^{\pi \Gamma_n}>0$, we have (on the event 
      $\{X_\infty=1\}$) $1-X_n\leq 
      Ce^{-\pi\Gamma_n}$ and $(1-X_n)/\g_{n+1}\leq 
      Ce^{-\pi\Gamma_{n+1}}/\g_{n+1}$. Then one can slightly relax the assumption in claim $(b)$ since it follows from 
      Lemma~\ref{LemmeFast} that if $\e_n\geq 0$ for $n$ large enough,
      $\sum_n(1-X_n)<+\infty$ almost 
      surely on $\{X_\infty=1\}$.
     \end{Rem}
 \noindent {\sc Proof of Lemma~\ref{LemmeFast}:} 
 Starting from (\ref{mean-algo-to1}), we have
 \begin{eqnarray}
     \frac{1-X_{n+1}}{\gamma_{n+1}}&=&
     \frac{1-X_n}{\gamma_{n+1}}-\pi X_n(1-X_n)-\D M_{n+1}\nonumber\\
     &=&(1-X_n)\left(\frac{1}{\gamma_n}+\e_n+\pi-\pi X_n\right)-\D M_{n+1}\nonumber\\
     &=&\frac{1-X_n}{\gamma_n}\left(1+\e_n\gamma_n+\pi 
     \gamma_n(1-X_n)\right)-\D M_{n+1}\label{star=},
 \end{eqnarray}
 so that, for $n\geq n_0$, $Z_{n+1}\geq Z_n-\D M_{n+1}$, which proves 
 that $(Z_n)_{n\geq n_0}$ is a submartingale.
 Now set
$\tau_{_L}:=\min\{n\ge n_0\,:\, 1-X_n>L\g_{n+1}\}$, $L>0$. 
Then the stopped submartingale $(Z^{\tau_{_L}}_n)_{n\ge n_0}$ satisfies
\[
(\Delta Z^{\tau_{_L}}_{n+1})_{_+}\le \mbox{\bf 1}_{\{\tau_{_L}\ge \,n+1\}}(\Delta Z_{n+1})_{_+}\le L +
\sup_n\|\Delta M_n\|_{\infty}.
\]
Consequently the sub-martingale $(Z^{\tau_{_L}}_n)_{n\ge n_0}$ is bounded  with bounded increments. Hence it 
converges ($\P_x$-$a.s.$ and in $L^1(\P_x)$) toward an integrable random variable $\zeta^L_{_\infty}$. Furthermore
(see~\cite{Neveu}) the conditional variance increment process of its martingale part also  converges  to a finite
random variable as  $n\to + \infty$. This reads
\[ 
\sum_{n=n_0+1}^{\tau_{_L}} \E((\Delta M_n)^2\,|\, {\cal F}_{n-1}) <+\infty\qquad \mbox{$\P_x$-$a.s.$}.
\]
But, we know from Proposition~\ref{ProVarCond} that
\[
 \E((\Delta M_n)^2\,|\, {\cal F}_{n-1}) \geq p_{_B} X_{n-1}(1-X_{n-1}).
\]
Consequently, 
\[
\{X_\infty = 1\}\cap \left(\cup_{p\in \N}\{\tau_{_p}=+\infty\}\right)\subset \{\sum_n 1-X_n <+\infty\}.
\] 
We conclude by observing  that $\cup_{p\in 
\N}\{\tau_{_p}=+\infty\}=\{\sup_n\frac{1-X_n}{\g_{n+1}}<+\infty\}.$\hfill$\cqfd$
 
\ms
 \noindent {\sc Proof of Proposition \ref{Fact3}:}  We first assume
 that $\sum_n\gamma_n\e_n^+<+\infty$. 
 The proof is based, as in Lemma~\ref{LemmeFast}, on the study of the sequence
 $((1-X_n)/\gamma_n)$. We deduce from (\ref{star=}) that
 \begin{eqnarray}
     \frac{1-X_{n+1}}{\gamma_{n+1}}
     &\leq &\frac{1-X_n}{\gamma_n}\left(1+\e_n^+\gamma_n+\pi 
     \gamma_n(1-X_n)\right)-\D M_{n+1}.\label{star}
 \end{eqnarray}
 Hence
 \begin{equation}
     \label{surmartingale}
     \E\left(\frac{1-X_{n+1}}{\gamma_{n+1}}
        \;|\;\F_n\right)\leq
        \frac{1-X_n}{\gamma_n}\left(1+\e_n^+\gamma_n+\pi 
     \gamma_n(1-X_n)\right).
 \end{equation}
 We know from Proposition~\ref{Fact1} that, on the set 
 $\{X_\infty=1\}$,
 we have $\displaystyle \sup_n(1-X_n)e^{\pi \Gamma_n}<+\infty$, so that
 $\gamma_n(1-X_n)\leq C\gamma_ne^{-\pi\Gamma_n}\leq C$ for some $C>0$,
 and $\sum_n \gamma_n(1-X_n)<+\infty$. We now deduce 
 from~(\ref{surmartingale}) and a supermartingale argument that,
 on $\{X_\infty=1\}$, the sequence $((1-X_n)/\gamma_n)_{n\in\N}$
 is almost surely convergent.
 
 On the other hand, with the notation $Z_n=(1-X_n)/\gamma_n$, we 
 know
 from (\ref{star}) that
 \[
 \D M_{n+1}\leq Z_n-Z_{n+1}+Z_n\left(\e_n^+\gamma_n+\pi 
     \gamma_n(1-X_n)\right).
 \]
 Therefore, on $\{X_\infty =1\}$ the martingale $M_n$ is bounded from
 above, and, since it has bounded jumps, we must have 
 $<\!M\!>_\infty<+\infty$ almost surely. We know from 
 Proposition~\ref{ProVarCond} that $<\!M\!>_\infty\geq p_{_B}\sum_n 
 X_{n-1}(1-X_{n-1})$. Hence $\sum_n(1-X_n)<+\infty$ a.s. on 
 $\{X_\infty=1\}$.
 
 \ms
 
 We now assume that $\liminf \e_n>0$, so that for $n$ large enough
 (say $n\geq n_0$), we have
 \begin{equation}\label{liminfpos}
 \frac{1}{\gamma_{n+1}}-\frac{1}{\gamma_n}-\pi\geq \e,
 \end{equation}
 for some $\e>0$. In particular the sequence $(\gamma_n)_{n\geq n_0}$
 is non-increasing and, for $n\geq n_0$,
 \[
 \gamma_n-\gamma_{n+1}\geq (\pi +\e)\gamma_n\gamma_{n+1},
 \]
 which implies $\sum_n\gamma_n^2<+\infty$.
 We also have, for $n\geq n_0$,
 \[
 \gamma_{n+1}\leq \gamma_n(1-(\pi+\e)\gamma_{n+1})\leq 
 e^{-(\pi +\e)\gamma_{n+1}}.
 \]
 Therefore, for $k\geq n\geq n_0$,
 \[
 \gamma_k\leq \gamma_ne^{-(\pi +\e)(\Gamma_k-\Gamma_n)},
 \]
 \begin{eqnarray*}
 \mbox{and }\hskip 3 cm     \sum_{k\geq n}\gamma_k^2 e^{\pi \Gamma_k}&\leq &
      \sum_{k\geq n}\gamma_k \gamma_n 
      e^{-(\pi +\e)(\Gamma_k-\Gamma_n)}e^{\pi \Gamma_k}\hskip 4 cm\\
      &=&\gamma_ne^{(\pi +\e)\Gamma_n}\sum_{k\geq 
      n}\gamma_ke^{-\e\Gamma_k}\\
      &\leq &\gamma_ne^{(\pi +\e)\Gamma_n}\int_{\Gamma_{n-1}}^\infty 
      e^{-\e x}dx\\
      &\leq &\gamma_n\frac{e^{\e\gamma_n}}{\e}e^{\pi \Gamma_n}.
 \end{eqnarray*}
We have thus proved not only that $\sum_n\gamma_n^2e^{\pi 
\Gamma_n}<+\infty$, but also that
\[
\sum_{k\geq n}\gamma_k^2 e^{\pi \Gamma_k}\leq C\gamma_ne^{\pi \Gamma_n}
\]
for some $C>0$. It then follows from Proposition~\ref{Fact2} that,
on the set $\{Y_\infty =0\}$, $(1-X_n)\leq 
C\theta_n\gamma_ne^{\pi\Gamma_n}$, and, using Remark~\ref{rem-theta},
we get $\displaystyle \sup_n(1-X_n)/\gamma_n<+\infty$ a.s. on $\{Y_\infty =0\}$.
We complete the proof by applying Lemma~\ref{LemmeFast}.
\hfill$\cqfd$
    
\begin{Rem}
        Assume, with the notation of Proposition~\ref{Fact3}, that
        $\liminf \e_n^+>0$ and $\sum_ne^{-p_{_A}\Gamma_n}<+\infty.$
        This is the case if $\gamma_n=C/(n+C)$, with $\pi C<1<p_{_A}C$. 
        Then, we deduce from Propositions~\ref{Fact3} and \ref{vitevitevers1}
        that $0<\P(Y_\infty=0)<1$ and that, on $\{Y_\infty=0\}$
        the sequence $(1-X_n)e^{p_{_A}\Gamma_n}$ converges to a positive 
        limit, whereas on $\{Y_\infty>0\}$, $(1-X_n)e^{\pi\Gamma_n}$ 
        converges to a positive limit almost surely.
    \end{Rem}
 \section{A parametric guide to the rates}\label{guidedtour}
In this section we will call {\em fast} a  rate of the algorithm which induces that the  error series  converges $i.e.$ 
$\sum_n 1-X_n<+\infty$ when $X_n\to 1$ and   
$\sum_n X_n<+\infty$  when $X_n\to 0$. Other rates will be considered as {\em slow}.

\ss
Assume (at least for large enough $n$)  that
\[
\g_n = \left(\frac{C}{C'+n}\right)^{\a},\quad \alpha\!\in(0,1],\quad C,\, C'>0.
\]
Then,  the algorithm behaves as follows: 
\begin{itemize}
\item If ($\alpha \!\in(0,1)$) or ($\alpha =1$ \& $Cp_{_B}>1$) then the algorithm is fallible with positive probability from any
$x\!\in[0,1)$ (note that  this probability is lower than $1$ if $x\!\in (0,1)$). When failing, it always goes to
$0$ at a fast rate, ($n^{-Cp_{_B}}$ if
$\alpha =1$). This follows from Proposition~\ref{EqFaillib}.
\item If $\alpha=1$ and $C\le \frac{1}{p_{_B}}$, the algorithm is infallible from any $x\!\in (0,1]$. This follows from
Proposition~\ref{ProConv0}$(b)$.
\end{itemize}
As concerns rates one has 
\begin{itemize}
\item If $\alpha=1$ and $C\ge\frac{1}{\pi}$ then the -- fast~-- rate of convergence is $n^{-Cp_{_A}}$ on $\{X_n\to 1\}$.
This follows from Proposition~\ref{Fact3}$(a)$.
\item If $\alpha=1$ and $\frac{1}{p_{_A}}< C<\frac{1}{\pi}$ then exactly {\em two rates of convergence  occur with positive
$\P_x$-probability} on
$\{X_n\to 1\}$: a slow one -- $n^{-C\pi}$~-- and a fast one -- $n^{-Cp_{_A}}$. This follows from
Proposition~\ref{vitevitevers1} and~\ref{Fact3}$(b)$ (see remark~7).
\item If   $\alpha=1$ and  $C\le \frac{1}{p_{_A}}$ then (the algorithm is infallible from any $x\!\in (0,1]$) but only the
slow rate  of convergence survives $i.e.$ $n^{-C\pi}$ on $\{X_n\to 1\}$. This follows from Proposition~\ref{vitevitevers1}.
\end{itemize}

\noindent Note as corollaries that, 

\ss
-- when $2\,p_{_B}\le p_{_A}$ (then $\frac{1}{\pi}\le \frac{1}{p_{_B}}$): it is possible to choose
$C\!\in[\frac{1}{\pi},\frac{1}{p_{_B}}]$  so that  the algorithm is simultaneously infallible and converging with
a fast rate. This is possible because in some sense $p_{_A}$ and $p_{_B}$ are remote enough.  The  fastest
achievable rate is  
$n^{-\frac{p_{_A}}{p_{_B}}}$  (with
$C =\frac{1}{p_{_B}}$).  Of course such  a specification is purely theoretical since
$p_{_A}$ and
$p_{_B}$ are supposed to be unknown.

\ss
-- when $\,p_{_B}<p_{_A}< 2\,p_{_B}$ (then $ \frac{1}{p_{_B}}<\frac{1}{\pi}$): there is no access to fast
converging rates within infallibility, because
$p_{_A}$ and
$p_{_B}$ are too close to each other .

\ss
-- in any case, when no information is available on the parameters $p_{_A}$
and $p_{_B}$, the ``blind" choice $C=1\le \frac{1}{p_{_A}}$  which ensures infallibility  induces a slow rate of convergence, namely
$n^{-\pi}$. In fact this rate can be very poor when $p_{_A}$ and $p_{_B}$ get close to each other.

\medskip
At this point the conclusion can be the following: the higher the parameter $C$ is, the faster the algorithm goes. But if $C$
is too high, it may go wrong. 
 
\medskip
-- One further point to be noticed is  that what we
called the slow rate -- $e^{-\pi\,\Gamma_n}$ -- for the algorithm is but {\em the rate of its mean deterministic version}
(see~\cite{LAPATA} for details). So, even when it is infallible (that is  converges to the same limit as its mean version),
it always converges  at least as fast as this deterministic procedure (which is of no practical interest since its implementation
would  require $p_{_A}$ and $p_{_B}$ to be known). When no information is available on the parameters
$p_{_A}$ and $p_{_B}$, this is the rate which is actually obtained.

\medskip As a conclusion, the   convergence rate behaviour of this stochastic approximation algorithm is completely non-standard.
Thus, from a mathematical viewpoint, one last feature to be noticed is the unusual ``spectrum" of the rates since the switching from
one rate to another  takes place ``progressively'' with a  range  of values of the parameter $C$ for the gain parameter  for which 
two different rates are achieved with positive probability.

\end{document}